\documentclass[12pt]{amsart}
\usepackage{amssymb,latexsym,amsmath,amsthm,amscd}
\usepackage{setspace}
\usepackage[all]{xy} \xyoption{arc}
\usepackage[left=3cm,top=2cm,right=3cm,bottom = 2cm]{geometry}
\usepackage{graphicx}
\usepackage[usenames,dvipsnames]{color}
\usepackage{charter}

\theoremstyle{plain}
\newtheorem{thm}{Theorem}

\theoremstyle{definition}

\theoremstyle{remark}
\newtheorem{rmk}[thm]{Remark}

\begin{document}
\title[Vertex superalgebras]{Most vertex superalgebras associated to an odd unimodular lattice of rank~$24$ have an $N{=}4$ superconformal structure}
\author{Gerald H\"{o}hn and Geoffrey Mason}
\address{Department of Mathematics, Kansas State University, Manhattan, KS 66506}
\email{gerald@monstrous-moonshine.de}
\address{Department of  Mathematics, University of California at Santa Cruz,Ca 95064}
\email{gem@ucsc.edu}

\thanks{The first author was supported by the Simons Foundation Collaboration Grant $\#355294$}
\thanks{The second author was  supported by the Simons Foundation $\#427007$}


\begin{abstract} Odd, positive-definite, integral, unimodular lattices $N$ of rank~$24$ were classified by Borcherds.\ 
There are $273$ isometry classes of such lattices.\ 
Associated to them are vertex superalgebras $V_N$ of central charge $c{=}24$.\ We show that 
\emph{at least\/} $267$ of these vertex operator superalgebras contain an $N{=}4$ superconformal subalgebra of central charge $c'{=}6$.\ 
This is achieved by studying embeddings $L^+{\subseteq} N$ of a certain rank~$6$ lattice $L^+$.
\end{abstract}
\maketitle
 

\section{Introduction}
The results discussed in this paper concern positive-definite, integral, unimodular lattices of rank $24$ and the vertex algebras that they define.\  
We  find it very convenient to call any such lattice a \emph{Niemeier} lattice, adding the adjective \emph{odd} or \emph{even} when appropriate.\ 
Thus for us an \emph{even Niemeier lattice} is what is now called simply a Niemeier lattice in the literature.\

\medskip
More often than not, the basic axiom for a vertex algebra, i.e., the \emph{Jacobi identity}, is supplemented by
 additional conditions such as the existence of a Virasoro field, or by a superconformal structure.\ In the case of vertex superalgebras there is a hierarchy of
 such structures, often denoted by $N{=}1$, $N{=}2$, $N{=}4$, $\ldots$ becoming more intricate with increasing $N$.\ 
For general background we refer the reader to \cite{DL}, \cite{LL}, \cite{K} and for an early comparison of these super structures, see \cite{GSW}.
 
 \medskip
Vertex algebras containing an $N{=}4$ superconformal structure of central charge~$6$, in particular, have come into recent prominence because of their intervention
in  \emph{Mathieu Moonshine} \cite{EOT}.\ Experience suggests that moonshine-like phenomena may be related to Niemeier lattices as we have defined them and
thus one is lead to consider odd Niemeier lattices $N$ of rank~$24$ and the vertex superalgebras
$V_N$ that they define.\
A natural question presents itself{:}\ does $V_N$ have an $N{=}4$ structure?\ 
In this paper we will show that most of the vertex operator superalgebras $V_N$  \emph{do} have an $N{=}4$ superconformal structure.
 
 \medskip
This question  was considered in \cite{MTY} for general lattice vertex  superalgebras $V_L$.\ 
However, the precise definition of the (small) $N{=}4$ algebra
is  awkward to handle  (cf.\ \cite{MTY}, Section~2) and makes the question hard to deal with.\ A way around this difficulty was proposed in \cite{MTY}.\
Namely, it was shown (\emph{loc. cit.}, Proposition~26) that a certain integral rank~$6$ lattice which we denote by $L^+$,
has the property  that $V_{L^+}$ contains an $N{=}4$ superconformal algebra $A$ and that $V_{L^+}$ and $A$ share the same Virasoro element 
(the canonical one defined by the Siegel-Sugawara construction).\ 
The lattice $L^+$ is spanned by elements $\alpha_1$, $\ldots$, $\alpha_6$ and $h {:=}\tfrac{1}{3}(\alpha_1{+}\cdots{+}\alpha_6)$
and equipped with a bilinear form satisfying $(\alpha_i, \alpha_j){=}3\,\delta_{ij}$,\ ($i$, $j=1$, $\ldots$,~$6$).

 
 \medskip
 In this way one can establish the existence of an $N{=}4$ superconformal subalgebra of central charge $6$ in $V_N$ by 
proving the existence of a sublattice  $L^+{\subseteq} N$, for we then have the tower $A{\subseteq}V_{L^+}{\subseteq}V_N$.\
 Our main results, obtained by a computer search, identify which odd  Niemeier lattices admit such an embedding.\ 

\section{Results}

 \medskip
In order to state our main results, we need to review the theory of Niemeier lattices~\cite{B,CS}.\
The important  idea of \emph{neighboring} lattices is due to Kneser.\ Suppose that $A$ and $B$ are a pair of equal rank lattices.\
We say that $A$ and $B$ are \emph{neighbors} in case their intersection $D$ (in a common Euclidean space $E$ containing both) has index~$2$ in each of them.\ 
If $A$ and $B$ are neighboring even Niemeier lattices then there is a unique odd Niemeier lattice $N$ that neighbors both $A$ and $B$.\ 
($A$, $B$ and $N$ are the lattices in $E$ containing $D$ with index~$2$).\
Furthermore, every odd Niemeier lattice $N$ arises in this way from an unordered pair $(A, B)$ of even Niemeier lattices that is unique up to isometry.\ We say that
$N$ has \emph{type} $(A, B)$.\ In general there will be several isometry classes of $N$ of the \emph{same type}.\ 
An even Niemeier lattice is characterized by its (semisimple) \emph{root system} (which by convention is empty in the case of the Leech lattice $\Lambda$).\ Thus 
the type of an odd Niemeier lattice $N$ is determined by the corresponding pair $(\Phi_A, \Phi_B)$ of root systems, which we also refer to as the \emph{type of $N$}.

 \medskip
We are also interested in the \emph{minimum norm}  $\mu$ of a Niemeier lattice $N$, i.e., the least value of
$(\alpha, \alpha)$ for $0{\not=}\alpha{\in}N$.\ If $N$ is even then $\mu{=}2$  with the single exception of the
Leech lattice $\Lambda$ which has $\mu{=}4$.\ For odd $N$ there are three possibilities as in the following table{:}

 \medskip
\begin{table}[htp]
\caption{Minimum norms for odd Niemeier lattices}
\begin{center}
\begin{tabular}{l|rrr}
$\mu$  & $1$ &  $2$  & $3$ \\ \hline
$\#$ of classes & $116$ & $156$ & $\phantom{22}1$
\end{tabular}
\end{center}
\label{default}
\end{table}

\medskip
The unique lattice of minimum norm~$3$ is called the \emph{odd Leech lattice} $\Lambda^{\rm odd}$; it has type $(\emptyset, A_1^{24})$.

\medskip
If $N$ is an odd Niemeier lattice of minimum norm~$1$, its two even Niemeier neighbors are necessarily isometric.\ 
This is because the Weyl reflection determined by a norm~$1$ vector \emph{exchanges} the two even neighbors.\ Thus the type of such a lattice always takes the form $(\Phi, \Phi)$ 
where $\Phi$ is the root system common to the even neighbors.
 
 \medskip
This situation dos not necessarily pertain in the case of odd Niemeier lattices  with $\mu{=2}$.\ 
The \emph{neighboring graph}, with nodes given by even Niemeier lattices and edges between $A$ and $B$ indexed by the odd lattices of type $(A, B)$
and with $\mu{=}2$ was first constructed by Borcherds~\cite{B}.\ It is Figure~17.1 in \cite{CS}.

\medskip
We can now state our main results.
\begin{thm}\label{thm1} Suppose that $N$ is an odd Niemeier lattice.\ Then the following hold:
\begin{enumerate} 
\item If $\mu{=}1$, then either $N$ is of type $(\emptyset, \emptyset)$, or else there is a sublattice $L^+\subseteq N$.
\item If $\mu{=}2$, exactly one of the following holds{:}
 \begin{eqnarray*}
 &&(a)\ \mbox{$N$ has type $(E_8^3, D_8^3)$, $(D_{16}E_8, A_{15}D_9)$, $(D_{12}^2, A_{12}^2)$ \ or $(D_{24}, A_{24})$},\\
 &&(b)\ \mbox{there exists a sublattice $L^+{\subseteq }N$}.
 \end{eqnarray*}
\item If $\mu{=}3$, the odd Leech lattice $\Lambda^{\rm odd}$ does not contain $L^+$.
\end{enumerate}
 \end{thm} 

\begin{rmk} (i)  There is a \emph{unique} odd Niemeier lattice for each of the types listed in Theorem \ref{thm1}(2)(a).\ 
 Thus there are exactly $152$ isometry classes of odd Niemeier lattices of minimum norm~$2$ that contain the lattice $L^+$.\\
(ii) The odd Leech lattice $\Lambda^{\rm odd}$ cannot possibly contain $L^+$  because the even sublattice of $\Lambda^{\rm odd}$ is contained in $\Lambda$ and hence has minimum norm ${\geq}4$, whereas $L^+$ contains vectors of norm~$2$ (namely ${\pm}h$).\\
(iii)\ For similar reasons, the unique lattice of type $(\emptyset, \emptyset)$ in Theorem \ref{thm1} (1) cannot contain $L^+$.\ Thus
we can paraphrase Theorem \ref{thm1} (1) as follows{:}\ every one of the $116$ odd Niemeier lattices  with $\mu{=}1$ that
\emph{could} contain $L^+$, \emph{does} contain $L^+$.
\end{rmk}

The proof of Theorem \ref{thm1} is computational.\ We used the neighborhood method as implemented in MAGMA~\cite{MAGMA} 
together with some optimizations to create the list of $273$ odd
Niemeier lattices together with their two even neighbors.\ We also computed the automorphism
group and checked the result by a mass formula (\cite{CS}, Chapter~16).

\medskip
The lattice $L^+$ has a basis consisting of the vector $h$ of norm~$2$ and the five pairwise orthogonal
norm~$3$ vectors $\alpha_1$, $\ldots$, $\alpha_5$ with scalar product $(\alpha_i,h)=1$.\
For each odd Niemeier lattice $N$, this permits  use of the following search method for the sublattice $L^+$:\
do a  backtracking search by selecting all possible  norm~$2$ vectors $h$ in $N$, then
for given vector $h$  take all norm~$3$ vectors $\alpha_1$ with $(h,\alpha_1)=1$, then for such pairs $\{ h,\alpha_1\}$ select all 
norm~$3$ vectors  $\alpha_2$ with $(h,\alpha_2)=1$
and $(\alpha_1,\alpha_2)=0$, then, in similar way, all possible vectors $\alpha_3$, $\alpha_4$ and $\alpha_5$.\ 
If we find a tuple $(h,\, \alpha_1,\, \ldots, \alpha_5)$ we stop our search for the lattice $N$.\
Otherwise we continue until all possibilities are exhausted.

\medskip
To speed up the computation, we used the automorphism group $O(N)$ of the lattice $N$.\
First we compute the orbits of $O(N)$ on the roots of norm~$2$.\ Then it is enough to
choose a random vector $h$ in each orbit.\ Then we compute the pointwise stabilizer $H_1$ of $h$
in $O(N)$ and decompose the set of norm~$3$ vectors $\alpha$ with $(h,\alpha)=1$ into $H_1$-orbits.\
Then it is enough to choose a random vector $\alpha_1$ in each such orbit.\ Then 
we compute the pointwise stabilizer $H_2$ of the two vectors $h$ and $\alpha_1$ in $O(N)$ and select a random vector
among the possible $\alpha_2$ in each $H_2$-orbit.\ For $\alpha_3$, $\alpha_4$ and $\alpha_5$ we went through  
all possible choices since this was quicker than computing additional stabilizers.

\medskip

\begin{rmk} 

(i) The computation took less than $24$ hours on a single processor machine.

(ii) We did not classify all orbits of sublattices $L^+\subseteq N$ under the action of $O(N)$, although this
can easily be achieved with our method.

(iii) We also tried a different computational approach.\ We constructed a single
embedding of $L^+$ into an odd Niemeier lattice $N$ and determined the orthogonal
complement lattice $K=(L^+\otimes {\bf R})^\perp \cap N$. Then we tried the neighborhood method
to compute all lattice $K'$ in the genus of $K$. By computing all gluings of all such lattices $K'$ with
$L^+$, we obtain odd Niemeier lattices $N'$ having $L^+$ as a sublattice and all such lattices can be constructed
that way.

\medskip
We found more than $10,000$ rank $18$ lattices $K'$ covering most of the genus of $K$. However, the
estimated run time for completing the calculation was more than a month.\ Apart from a few missing cases, 
the list of lattices $N'$ obtained in this way was identical with the one from our theorem.    
\end{rmk}

The following two questions remain open:

(i) Do the four exceptional lattice vertex superalgebras $V_N$  in part (2) of Theorem \ref{thm1} have an $N{=}4$ super conformal structure not
coming from a lattice embedding $L^+\subseteq N$?

(ii) Do most of the self-dual lattice vertex superalgebras of central charge~$24$ (cf.~\cite{H}),
which are not of lattice type have an $N{=}4$ superconformal structure?


\end{document}